\documentclass[12pt]{article}
\usepackage{amsmath}
\usepackage{amssymb}
\usepackage{mathtools}
\usepackage{graphicx}
\usepackage{hyperref}
\newtheorem{Theorem}{Theorem}
\newtheorem{Lemma}{Lemma}
\newtheorem{Corollary}{Corollary}

\newtheorem{Remark}{Remark}
\newcommand{\dd}{{\rm d}}
\newcommand{\ee}{{\rm e}}
\newcommand{\eps}{\varepsilon}

\begin{document}

\begin{center}
{ \Large On the Discontinuous Galerkin Finite Element Method}\vspace{2mm}

{\Large for Reaction--Diffusion Problems:}\vspace{2mm}

{\Large Error Estimates in Energy and Balanced Norms}\vspace{10mm}
  
{\large Helena Zarin\footnote{Department of Mathematics and Informatics, Faculty of Sciences, University of Novi Sad, 21000 Novi Sad, Serbia, helena.zarin@dmi.uns.ac.rs}, 
Hans--G\"org Roos\footnote{Institut f\"ur Numerische Mathematik, Technische Universit\"at Dresden, 01062 Dresden, Germany, hans-goerg.roos@tu-dresden.de}
}\vspace{10mm}
  
\end{center}

\begin{abstract}
A nonsymmetric discontinuous Galerkin FEM with interior penalties has been applied to one--dimensional singularly perturbed reaction--diffusion problems. Using higher order splines on Shishkin--type layer--adapted meshes and certain graded meshes, robust convergence has been proved in the corresponding energy norm and in a balanced norm. Numerical experiments support theoretical findings.\vspace{3mm}

finite element method, layer--adapted mesh, balanced norm

\end{abstract}\vspace{2mm}

%\end{keyword}

%\end{frontmatter}

\section{Introduction}

Singularly perturbed problems have been extensively studied over the last few decades. In a vast literature, different numerical methods for constructing robust numerical approxi\-mations have been presented; see e.g. the books \cite{L10,M02,MOS96,RST08}, survey papers \cite{R12,S05}, and references therein. In the context of finite element methods, error bounds have been usually derived in energy norms associated to corresponding bilinear forms. However, in several recent papers \cite{FR13,LS12,MX16,RS14}, the weakness of the energy norm to recognize characteristics of the layers has been addressed. Thus, a stronger {\em balanced} norm has been introduced, in which both regular and layer solution components are uniformly bounded. 

Here we are interested in numerical solving of a reaction--diffusion problem using a nonsymmetric version of the discontinuous Galerkin finite element method with interior penalties (NIPG method, \cite{HSS02,RZ03}). The purpose of the paper is twofold. First, we prove a parameter--uniform convergence in an energy norm extending the analysis from \cite{RZ03} to higher order splines on a class of Shishkin--type meshes and graded meshes of Duran--Lombardi type. Second and more important, we also prove error estimates in a balanced norm. For reaction--diffusion problems, so far those error estimates exist only for the Galerkin finite element method \cite{RS14}, certain mixed methods \cite{LS12}, and an $hp$ finite element method on a spectral boundary layer mesh \cite{MX16}. 

In order to clearly present basic ideas, here we focus on a one--dimensional reaction--diffusion problem, while at the end of the paper we address the questions of genera\-lizations to systems of reaction--diffusion equations and the two--dimensional case. Our model problem is the following singularly perturbed differential equation
\begin{equation}\label{eq:dspp}
\left\{
\begin{array}{c}
-\eps^2u''(x)+c(x)u(x)=f(x), \quad x\in\Omega=(0,1), \\[1ex]
u(0)=u(1)=0,
\end{array}
\right.
\end{equation}
where $0<\eps\ll1$ is a perturbation parameter, $c,f$ are smooth functions  such that
\begin{equation}\label{eq:ass}
c(x)\geq \tilde\gamma^2>0
	\qquad \mbox{for all\,\,}x\in\overline\Omega:=[0,1],
\end{equation}
with some constant $\tilde\gamma$. The behaviour of the solution $u$ to \eqref{eq:dspp}--\eqref{eq:ass} and its derivatives is already known, \cite{L10}: the solution has two boundary layers and there exists a solution decomposition $u=S+E$, where
\begin{equation}\label{eq:decomp}
|S^{(k)}(x)| \leq C, \qquad 
|E^{(k)}(x)| \leq C\eps^{-k} \left(\ee^{-x\gamma/\eps}+\ee^{-(1-x)\gamma/\eps}\right),
\end{equation}
for all $x\in\Omega$, $\gamma\in(0,\tilde\gamma)$, and $k=0,1,\ldots,q$ (the order $q$ depends on the smoothness of the data). Beside in error analysis, this a priori information on the solution influences the construction of a discretization mesh that should adequately resolve the layers.

The paper is organized as follows. In the next section we describe the NIPG method as well as layer--adapted meshes of Shishkin--type (S--type) and recursively defined graded meshes. In Section~\ref{sec:error} we present the analysis of the method in both energy and balanced norms, separately estimating interpolation and discretization errors. Section~\ref{sec:num} contains the results of numerical experiments in order to illustrate theoretical bounds. In the last section we comment on more general problems and analysis extensions.

{\em Notation.} With $C$ we denote a generic positive constant independent of the perturbation parameter $\eps$ and the number of degrees of freedom $N$. For a set $D\subseteq\Omega$, $\overline D$ denotes its closure, and $\mathbb{P}_k(D)$ is the set of polynomials defined on $D$, of the highest degree $k\ge 1$. Moreover, on $D$ we use the standard notation for Banach spaces $L^q(D)$, Sobolev spaces $H^q(D)$, norms $\|\cdot\|_{L^q(D)}$, $\|\cdot\|_{H^q(D)}$ and seminorm $|\cdot|_{H^q(D)}$. The scalar product in $L^2(D)$ is denoted with $(\cdot,\cdot)_D$; we write $(\cdot,\cdot)$ when $D=\Omega$.

\section{Problem discretization}

\subsection{The NIPG method}

Let $N\geq 4$ be an even integer, and $\{x_0,x_1,\ldots,x_{N}\}$ a general mesh on $\overline\Omega$ that defines elements $I_i=(x_{i-1},x_i)$ such that $\overline\Omega=\bigcup_{i=1}^{N}\overline{I_i}$. We take $x_0=0$, $x_N=1$. Our broken Sobolev space will be $V=\left\{v\in L^2(\Omega):\left.v\right|_{I_i}\in H^1(I_i),\,i=1,2,\ldots,N\right\}$. For a function  $v\in V$, a jump $[v]_i$ and an average $\langle v\rangle_i$ at the mesh node $x_i$ are defined with $[v]_i=v(x_i+0)-v(x_i-0)$, $\langle v\rangle_i=\left(v(x_i+0)+v(x_i-0)\right)/2$, $i=1,2,\ldots,N-1$, $[v]_0=\langle v\rangle_0=v(x_0+0)$, $[v]_{N}=-\langle v\rangle_{N}=-v(x_{N}-0)$.

Now, the weak formulation related to the NIPG method for the problem \eqref{eq:dspp}--\eqref{eq:ass} is: find $u\in V$ such that
\begin{equation}\label{eq:weak_form}
a(u,v)=l(v), \qquad \mbox{for all}\quad  v\in V,
\end{equation}
where for functions $w,v\in V$,
\begin{align}
a(w,v)	 & =a_G(w,v)+\sum_{i=0}^{N}
	 \left(\eps^2\langle w'\rangle_i[v]_i-\eps^2[w]_i\langle v'\rangle_i+\sigma_i[w]_i[v]_i\right), \label{eq:blf} \\
a_G(w,v) & = \sum_{i=1}^{N}\int_{I_i}
	\left(\eps^2w'(x)v'(x)+c(x)w(x)v(x)\right)\dd x,
		 \nonumber  \\
l(v) & =\sum_{i=1}^{N}\int_{I_i}f(x)v(x)\dd x,
	 \nonumber
\end{align}
with penalty parameters $\sigma_i>0$ as constants for controlling the jumps of a discrete solution.

If the finite element space $V^N\subset V$ consists of $k$th degree piecewise polynomials defined on our general mesh, i.e. $V^N=\left\{v\in L^2(\Omega):\left.v\right|_{I_i}\in \mathbb{P}_k(I_i),\,i=1,2,\ldots,N\right\}$, then the discrete problem reads: find $u^N\in V^N$ such that
\begin{equation}\label{eq:discr_form}
a(u^N,v^N)=l(v^N), \qquad \mbox{for all}\quad  v^N\in V^N.
\end{equation}
Both \eqref{eq:weak_form} and \eqref{eq:discr_form} admit a unique solution due to the assumption \eqref{eq:ass}, and the bilinear form \eqref{eq:blf} defines an energy norm $\|w\|_{dG}^2:=a(w,w)$.

If $u^*\in V^{N}$ represents some interpolant (respectively projection) of the solution $u$, the analysis of $\|u-u^N\|_{dG}$ will emanate from the triangle inequality applied to the error decomposition
\begin{equation}\label{eq:err_decomp}
u-u^N=\eta+\chi,\qquad \eta:=u-u^*, \qquad \chi:=u^*-u^N.
\end{equation}
In the sequel we will use as well the standard (globally continuous) Lagrange interpolant $u^I$ as the (discontinuous) piecewise
$L^2-$projection $u^\pi$.

\subsection{Layer--adapted meshes}\label{ssec:mesh}

We consider as well Shishkin--type meshes, \cite{L10,RL99,RST08}, as graded meshes due to Duran and Lombardi, \cite{DL06}. Remark that it would also be possible to handle modified S--type meshes in the sense of \cite{FX16}, which include exponentially graded meshes from \cite{CX15}.

\subsubsection{S--type meshes}
For the given integer $N\geq 4$ divisible by 4, let $\lambda=(k+1)\eps/\gamma\ln N<1/4$,
be a mesh transition parameter of Shishkin type and let us assume $\eps\leq CN^{-1}$. Notice that the layer component $E$ of the solution in \eqref{eq:decomp} has the property 
\[
\max\left\{|E^{(k)}(\lambda)|,|E^{(k)}(1-\lambda)|\right\}\leq C\eps^{-k}N^{-(k+1)}.
\]

In order to adequately resolve the layers of the solution of \eqref{eq:dspp}, we construct the mesh such that it is equidistant on $\overline\Omega_c$ with the mesh step size $2(1-2\lambda)N^{-1}$ and gradually divided on $\overline\Omega_f$, where
\[
\Omega_c  = (\lambda,1-\lambda), \qquad 
\Omega_f = (0,\lambda)\cup(1-\lambda,1). 
\]
We choose transition points as $x_{N/4}=\lambda$, $x_{3N/4}=1-\lambda$. Following \cite{L10}, the layer--adapted mesh on $\overline\Omega_f$ is defined using two mesh generating functions $\phi_{1,2}$ that are continuous, piecewise continuously differentiable, $\phi_1$ ($\phi_2$) is monotonically increasing (decreasing), $\phi_1(0)=\phi_2(1)=0$, $\phi_1(1/4)=\phi_2(3/4)=\ln N$. Finally, the mesh points are
\begin{equation}\label{eq:mesh}
x_i=\begin{cases}
(k+1)\frac{\eps}{\gamma}\phi_1(i N^{-1}), & i=0,1,\ldots ,N/4, \\[1ex]
\lambda+2(1-2\lambda)(i N^{-1}-\frac{1}{4}), & i=N/4+1,\ldots,3N/4, \\[1ex]
1-(k+1)\frac{\eps}{\gamma}\phi_2(i N^{-1}), & i=3N/4+1,\ldots ,N.
\end{cases}
\end{equation}

In the sequel we assume the mesh generating functions satisfy
\begin{equation}\label{eq:mesh_assump}
N^{-1}\max|\phi'|\leq C, \quad \min_{k}\left(\phi\Big(\frac{k}{N}\Big)-\phi\Big(\frac{k-1}{N}\Big)\right)\geq CN^{-1},
\end{equation}
and define mesh characterizing functions $\psi=\ee^{-\phi}$ omitting indices for the mere of simplicity. In Table~\ref{tbl:mesh} we present examples of $\psi$ for different layer--adapted meshes from \cite{L10,RL99}: Shishkin mesh (S--), polynomial Shishkin mesh with $m>1$ (pS--), Bakhvalov--Shishkin mesh (BS--), and modified Bakhvalov--Shishkin mesh with $q=1/2+1/(2\ln N)$ (mBS--mesh).

\begin{table}[t!]
\begin{center}
\caption{Mesh characterising functions $\psi=\ee^{-\phi}$}\vspace{2mm}
{\renewcommand{\arraystretch}{1.1}
\begin{tabular}{|c|l|l|l|l|} \hline
& S--mesh & pS--mesh & BS--mesh &  mBS--mesh \\ \hline 
$\psi_1(t)$ & $N^{-4t}$ & $N^{-(4t)^m}$ & $1-4(1-N^{-1})t$ & $\ee^{-2t/(q-2t)}$ \\[1mm]
$\psi_2(t)$ & $N^{-4(1-t)}$ & $N^{-(4(1-t))^m}$ & $1-4(1-N^{-1})(1-t)$ & $\ee^{-2(1-t)/(q-2+2t)}$ \\[1mm]\hline
$\max|\psi'|$ & $C\ln N$ & $C(\ln N)^{1/m}$ & $C$ & $C$\\ \hline
$\max|\phi'|$ & $C\ln N$ & $Cm\ln N$ & $CN$ & $C\ln^2N$\\ \hline
\end{tabular} \label{tbl:mesh}
}
\end{center}
\end{table}

Following the technique from \cite{RL99}, the mesh step sizes $h_i=x_i-x_{i-1}$ satisfy
\begin{equation}\label{eq:h_i}
\begin{array}{rcll}
C(k+1)\eps N^{-1} \hspace{-2mm}& \leq h_i \leq  &\hspace{-2mm} 
	C(k+1)\eps N^{-1}\max|\phi'|, & \quad \mbox{on\,\,\,}\Omega_f, \\[1ex]
N^{-1} \hspace{-2mm}& \leq h_{i} \leq &\hspace{-2mm} 2N^{-1}, & \quad \mbox{on\,\,\,}\Omega_c.
\end{array}
\end{equation}

\subsubsection{Graded meshes}

Recursively generated meshes appear relatively often in the literature. In 1D, recursively gene\-rated meshes for a problem with a boundary layer characterized by the parameter $\eps$ and a layer width of order  $\eps$, have the form
\begin{align*}
 x_0&=0, \\
 x_1&=\eps H, \\
 x_i&=x_{i-1}+g(\eps, H, x_{i-1}), \quad i=2,\dots,M,
\end{align*}
with some parameter $H\in(0,1)$.

Following a proposal of Duran and Lombardi \cite{DL06}, we take the simplest mesh of that type
\begin{equation} \label{eq:graded}
x_i=\begin{cases}
iH\eps,  & \mbox{if}\quad 0\le i< [\frac{1}{H}]+1, \\[1ex]
(1+H)x_{i-1},  & \mbox{if}\quad [\frac{1}{H}]+1\le i\le M-1,\\[1ex]
1/2, & \mbox{if}\quad i=M,
\end{cases}
\end{equation}
where $M$ is such that 
\begin{equation}\label{eq:cond}
x_{M-1}<1/2 \quad \mbox{ and } \quad (1+H)x_{M-1}\ge 1/2. 
\end{equation}
In $[1/2,1]$ we use the same (reflected) mesh, i.e. $x_{M+i}=1-x_{M-i}$, $i=1,\dots,M$. The total number of mesh subintervals is $N=2M$. In case the last interval $(x_{M-1},1/2)$ is too small compared to $(x_{M-2},x_{M-1})$, we simply omit the mesh points $x_{M-1},x_{M+1}$, cf. \cite{DL06}.

Let $\ell=[\frac{1}{H}]$. The mesh step sizes $h_i=x_i-x_{i-1}$, $i=1,2,\dots,N$, satisfy
\begin{alignat*}{3}
h_{M+i} & =h_{M-i+1}, & \qquad & 1\le i\le M, \\
h_i & =H\eps, & \qquad & 1\le i\le \ell \quad \mbox{and} \quad N-\ell+1\le i\le N, \\
h_i & \le Hx_{i-1}\le Hx,  & \qquad & \ell+1\le i\le M, \quad x\in\overline{I}_i, \\
h_i & \le H(1-x_i)\le H(1-x), & \qquad & M+1\le i\le N-\ell, \quad x\in\overline{I}_i.
\end{alignat*}
These properties can be easily derived; e.g. the last inequality for indices $i=M+j$, $1\le j\le M-\ell$, follows from
\[
h_i=h_{M-j+1}\le Hx_{M-j}=H(1-x_{M+j})=H(1-x_i)\le H(1-x), \quad x\in\overline{I}_i.
\]
Moreover, the mesh step sizes can be estimated with  $CH\eps\le h_i\le H$, $i=1,2,\dots,N$.  

\begin{Remark}\label{rem}
On recursively generated meshes, the number of degrees of freedom $N$ is not known in advance. On the Duran--Lombardi mesh it is determined from \eqref{eq:cond}, and together with the perturbation parameter $\eps$ and the mesh parameter $H$, satisfies 
\[
    H\le C\,\frac{\ln(1/\eps)}{N}.
\]
\end{Remark}

\section{Error analysis}\label{sec:error}

\subsection{The interpolation error}

Similarly to \cite{RL99}, we can prove the following assertion on different norms of the interpolation error, considering some of them
elementwise. If the norm has to be considered elementwise, we characterize it by some index $d$.

\begin{Lemma}\label{thm:int_error}
For the projection error $\eta=u-u^*$ between the solution $u$ of the problem \eqref{eq:dspp}--\eqref{eq:ass} and its interpolant $u^*\in V^N$, on the Shishkin--type mesh \eqref{eq:mesh} it holds
\begin{alignat*}{4}
\|\eta\|_{L^\infty(\Omega_c)} & \leq CN^{-(k+1)}, &
	\quad \|\eta\|_{L^\infty(\Omega_f)} & \leq C\left(N^{-1}\max|\psi'|\right)^{k+1}, \\
\eps\|\eta'\|_{L_d^\infty(\Omega_c)} & \leq CN^{-k}, &
	\quad \eps\|\eta'\|_{L_d^\infty(\Omega_f)} & \leq C\left(N^{-1}\max|\psi'|\right)^k,  \\
\|\eta\|_{L^2(\Omega_c)} & \leq CN^{-(k+1)}, &
	\quad \|\eta\|_{L^2(\Omega_f)} & \leq C\eps^{1/2}\left(N^{-1}\max|\psi'|\right)^{k+1}, \\
\eps^{1/2}|\eta|_{H_d^1(\Omega_c)} & \leq CN^{-k},  &
	\quad \eps^{1/2}|\eta|_{H_d^1(\Omega_f)} & \leq C\left(N^{-1}\max|\psi'|\right)^k.
\end{alignat*}
\end{Lemma}

Let us choose interior penalty parameters $\sigma_i$, $i=0,1,\ldots,N$, as
\begin{equation}\label{eq:sigma}
\sigma_i=\begin{cases}
\eps,  & \mbox{if}\quad x_i\in\{x_0,x_{N}\}, \\
\eps N,  & \mbox{if}\quad x_i\in\Omega_c, \\
\eps N \left(\max|\psi'|\right)^{-1}, & \mbox{if}\quad x_i\in\Omega_f^\ast:=\overline\Omega_f\setminus\{x_0,x_{N}\}.
\end{cases}
\end{equation}
From the previous Lemma it follows
\begin{equation}\label{eq:int_error}
\|\eta\|_{dG} \leq CN^{-(k+1)}+C\eps^{1/2}N^{-k}\max|\psi'|^{k+1/2}.
\end{equation}

\begin{Remark}
For a globally continuous Lagrange interpolant $u^I\in V^N$ that satisfies $u^I(x_i\pm 0)=u(x_i\pm 0)$, we easily get 
\[
\|\eta\|_{dG}\leq CN^{-(k+1)}+C\eps^{1/2} \left(N^{-1}\max|\psi'|\right)^k,
\]
without making any specific choice for the penalization parameters.
\end{Remark}\smallskip

Next we consider the projection error on the graded mesh \eqref{eq:graded} that can be bounded similarly to \cite{BZ16,DL06}.

\begin{Lemma}\label{thm:int_error2}
For the projection error $\eta=u-u^*$ between the solution $u$ of the problem \eqref{eq:dspp}--\eqref{eq:ass} and its interpolant $u^*\in V^N$, on the Duran--Lombardi mesh \eqref{eq:graded} it holds
\begin{alignat*}{4}
\max\left\{\|\eta\|_{L^\infty(\Omega)},\|\eta\|_{L^2(\Omega)}\right\} & \leq CH^{k+1}, \\
\max\left\{\eps\|\eta'\|_{L_d^\infty(\Omega)},\eps^{1/2}|\eta|_{H_d^1(\Omega)}\right\} & \leq CH^k.
\end{alignat*}
\end{Lemma}\smallskip

Choosing the penalty parameters as
\begin{equation}\label{eq:sigma2}
\sigma_i=\begin{cases}
\eps,  & \mbox{if}\quad x_i\in\{x_0,x_{N}\}, \\
\eps H^{-1},  & \mbox{if}\quad x_i\in\Omega, \\
\end{cases}
\end{equation}
we get in the dG--norm
\begin{equation}\label{eq:int_error2}
\|\eta\|_{dG} \leq CH^{k+1}+C\eps^{1/2}N^{1/2}H^{k+1/2} \le CH^{k+1}+C\left(\eps\ln(1/\eps)\right)^{1/2}H^k.
\end{equation}
Here we have used the relation between $\eps$, $H$ and $N$ from Remark~\ref{rem}, cf. \cite{DL06}.

\subsection{The discretization error}

Next we estimate $\chi:=u^*-u^N$.
The Galerkin orthogonality property of the bilinear form \eqref{eq:blf} leads to $\|\chi\|_{dG}^2=-a(\eta,\chi)$. In the sequel we estimate each of the terms participating in $|a(\eta,\chi)|$, cf.~\eqref{eq:blf}.

We start with Shishkin--type meshes. Similarly to \cite{RZ03}, the Cauchy--Schwarz inequality and Lemma \ref{thm:int_error} imply
\begin{equation}
|a_G(\eta,\chi)| \leq C\left(\eps^{1/2}\left(N^{-1}\max|\psi'|\right)^k+N^{-(k+1)}\right)\|\chi\|_{dG}. \label{eq:12}
\end{equation}
When the interior penalty parameters are chosen as in \eqref{eq:sigma}, then Lemma \ref{thm:int_error} yields
\begin{equation}\label{eq:3}
\left|\eps^2\sum_{i=0}^{N}\langle \eta'\rangle_i[\chi]_i\right|
	\leq \left(\sum_{i=0}^{N}\eps^4\sigma_i^{-1}\langle \eta'\rangle_i^2\right)^{1/2}\|\chi\|_{dG}
	\leq C\eps^{1/2}N^{-k}\max|\psi'|^{k+1/2}\|\chi\|_{dG}.
\end{equation}
From the properties of $h_i$ from \eqref{eq:h_i} and interpolation error estimates we get
\begin{align*}
\sum_{i=1}^{N}\frac{\eps^2}{h_{i}}[\eta]_i^2
  & \leq C\sum_{x_i\in\Omega_c}\frac{\eps^2}{h_{i}}N^{-2(k+1)}
 	+C\sum_{x_i\in\Omega^\ast_f\cup\{x_{2N}\}}\frac{\eps^2}{h_{i}}\left(N^{-1}\max|\psi'|\right)^{2(k+1)} \\
  & \leq C\eps N^{-2k}\max|\psi'|^{2(k+1)}.
\end{align*}
Now, using inverse inequalities for the function $\chi\in V^N$, we derive
\begin{align}\label{eq:4}
\left| \eps^2\sum_{i=0}^N[\eta]_i\langle \chi'\rangle_i\right|
 &  \leq C\left(\left(\sum_{i=1}^N\frac{\eps^2}{h_{i}}[\eta]_i^2\right)^{1/2}+
 \left(\sum_{i=0}^{N-1}\frac{\eps^2}{h_{i+1}}[\eta]_i^2\right)^{1/2}\right)\|\chi\|_{dG} \nonumber \\[1ex]
 & \leq C\eps^{1/2}N^{-k}\max|\psi'|^{k+1}\|\chi\|_{dG}.
\end{align}
Finally,
\begin{equation}\label{eq:5}
\left|\sum_{i=0}^N\sigma_i[\eta]_i[\chi]_i\right|
	\leq \left(\sum_{i=0}^N\sigma_i[\eta]_i^2\right)^{1/2}\|\chi\|_{dG}
	\leq C\eps^{1/2}N^{-k}\max|\psi'|^{k+1/2}\|\chi\|_{dG}.
\end{equation}
Collecting \eqref{eq:12}--\eqref{eq:5} into $|a(\eta,\chi)|$, we obtain the estimate for the discretization error
\begin{equation}\label{eq:discr_error}
\|\chi\|_{dG}\leq CN^{-(k+1)}+C\eps^{1/2}N^{-k}\max|\psi'|^{k+1}.
\end{equation}

\begin{Remark}
When the interpolant is continuous, i.e., $[\eta]_i=0$, $i=0,1,\ldots,N$, previous ana\-lysis simplifies to a great extent. While \eqref{eq:12} remains, taking e.g. $\sigma_i=N$, $i=0,1,\ldots,N$, \eqref{eq:3} can be estimated with $C\eps\left(N^{-1}\max|\psi'|\right)^k\|\chi\|_{dG}$, implying 
\[
\|\chi\|_{dG}\leq CN^{-(k+1)}+C\eps^{1/2}\left(N^{-1}\max|\psi'|\right)^k.
\]
\end{Remark}\smallskip

The main result on the $\eps-$uniform convergence of the NIPG method \eqref{eq:discr_form} in the
energy norm on the discretization  mesh \eqref{eq:mesh} immediately follows from \eqref{eq:err_decomp}, \eqref{eq:int_error} and \eqref{eq:discr_error}.

\begin{Theorem}\label{thm:main}
Let  $u$ be the solution of the problem \eqref{eq:dspp}--\eqref{eq:ass} and $u^N\in V^N$ its numerical approximation that solves the discrete problem \eqref{eq:discr_form} on the layer--adapted mesh \eqref{eq:mesh}. If the penalty parameters are chosen as in \eqref{eq:sigma}, then
\[
 \|u-u^N\|_{dG}\leq CN^{-(k+1)}+C\eps^{1/2}N^{-k}\max|\psi'|^{k+1}.
 \]
\end{Theorem}\smallskip

As previously mentioned, the energy norm appears to be inadequate for detecting layer effects. For example, in our case one of the layer functions $\ee^{-x\gamma/\eps}$ from \eqref{eq:decomp} has the property $|\ee^{-x\gamma/\eps}|_{H^1(\Omega)}\leq C\eps^{-1/2}$. If in the energy norm we replace $\eps|\cdot|_{H_d^1(\Omega)}$ with $\eps^{1/2}|\cdot|_{H_d^1(\Omega)}$, we obtain the so--called balanced norm $\|\cdot\|_{dG,b}$.

Considering \eqref{eq:12}--\eqref{eq:5}, we observe that the term $N^{-(k+1)}$ without the factor $\eps^{1/2}$ arises from the estimate of $(c\eta,\chi)$ in the Galerkin part. But if we choose $u^*$ to be the (generalized) $L^2-$projection, defined by
\[
   (cu^{\pi},\xi)=(cu,\xi), \qquad \mbox{for all}\quad  \xi\in V^N,
\]
that term disappears and we have
\[
    \eps|\chi|_{H_d^1(\Omega)}\le \|\chi\|_{dG}\le C\eps^{1/2}N^{-k}\max|\psi'|^{k+1}.
\]
Consequently, we immediately get an error estimate in the balanced norm.

\begin{Corollary}\label{cor:B}
Let  $u$ be the solution of the problem \eqref{eq:dspp}--\eqref{eq:ass} and $u^N\in V^N$ its numerical approximation that solves the discrete problem \eqref{eq:discr_form} on the layer--adapted mesh \eqref{eq:mesh}. If the penalty parameters are chosen as in \eqref{eq:sigma}, then
\[
 \|u-u^N\|_{dG,b}\leq CN^{-k}\max|\psi'|^{k+1}.
 \]
\end{Corollary}\smallskip

Next we consider the error estimation on a graded mesh. Here we expect a weaker result (a weak dependence on $\eps$ in
the final error estimate). But the graded mesh of Duran--Lombardi has its advantages: it is not necessary to define
some transition point of the mesh; moreover, the mesh is robust in the sense that a mesh generated for a certain value
of $\eps$ can also be used for larger values of $\eps$.

We proceed in the same way as on a Shishkin--type mesh. First the Galerkin part yields
\[
|a_G(\eta,\chi)| \leq C\left(\eps^{1/2}H^k+H^{k+1}\right)\|\chi\|_{dG}. 
\]
For the terms corresponding to \eqref{eq:3}--\eqref{eq:5} we get bounds of the structure ${\cal{O}}(\eps^{1/2}N^{1/2}H^{k+1/2})$.
Following \eqref{eq:int_error2}, consequently we obtain

\begin{Theorem}\label{thm:main2}
Let  $u$ be the solution of the problem \eqref{eq:dspp}--\eqref{eq:ass} and $u^N\in V^N$ its numerical approximation that solves the discrete problem \eqref{eq:discr_form} on the Duran--Lombardi mesh \eqref{eq:graded}. If the penalty parameters are chosen as in \eqref{eq:sigma2}, then
\[
\|u-u^N\|_{dG}\leq CH^{k+1}+C\left(\eps\ln(1/\eps)\right)^{1/2}H^k.
\]
\end{Theorem}

The previous result can be restated in terms of the mesh node numbers $N$. Thus, employing Remark \ref{rem} the previous inequality reads
\begin{equation}\label{eq:main2}
\|u-u^N\|_{dG}\leq CH^k\le CN^{-k}\left(\ln(1/\eps)\right)^k.
\end{equation}
Clearly, the logarithmic dependence of the mesh parameter $H$ deteriorates the order of convergence as the polynomial degree increases.

Analogously as above we can also estimate the error in a balanced norm choosing the $L^2-$projection
as interpolant.

\begin{Corollary}\label{cor:2}
Let  $u$ be the solution of the problem \eqref{eq:dspp}--\eqref{eq:ass} and $u^N\in V^N$ its numerical approximation that solves the discrete problem on the Duran--Lombardi mesh \eqref{eq:graded}. If the penalty parameters are chosen as in \eqref{eq:sigma2}, then
\[
\|u-u^N\|_{dG,b} \leq  CH^{k}\left(\ln(1/\eps)\right)^{1/2}\leq CN^{-k}\left(\ln(1/\eps)\right)^{k+1/2}.
\]
\end{Corollary}

\section{Numerical results}\label{sec:num}

In this section we present the results of numerical experiments for the NIPG method \eqref{eq:discr_form} applied to layer--adapted meshes from Subsection~\ref{ssec:mesh}. The test problem is
\begin{equation}\label{eq:spp-test}
\left\{
\begin{array}{c}
-\eps^2u''(x)+(3-x^2)u(x)=f(x), \quad x\in(0,1), \\[1ex]
u(0)=u(1)=0,
\end{array}
\right.
\end{equation}
with the function $f$ such that \eqref{eq:spp-test} has the exact solution
\[
u(x)=\frac{\ee^{-x/\eps}+\ee^{-(1-x)/\eps}}{1+\ee^{-1/\eps}}-1+x^2(1-x)^2.
\]

In all our experiments we take the perturbation parameter $\eps=2^{-20}$. This choice is sufficiently small to bring out the singularly perturbed nature of \eqref{eq:spp-test}. Moreover, all integrals are approximated with the $5-$point Gauss--Legendre quadrature.

\begin{figure}[h!]
\begin{center}
\includegraphics[width=0.95\textwidth]{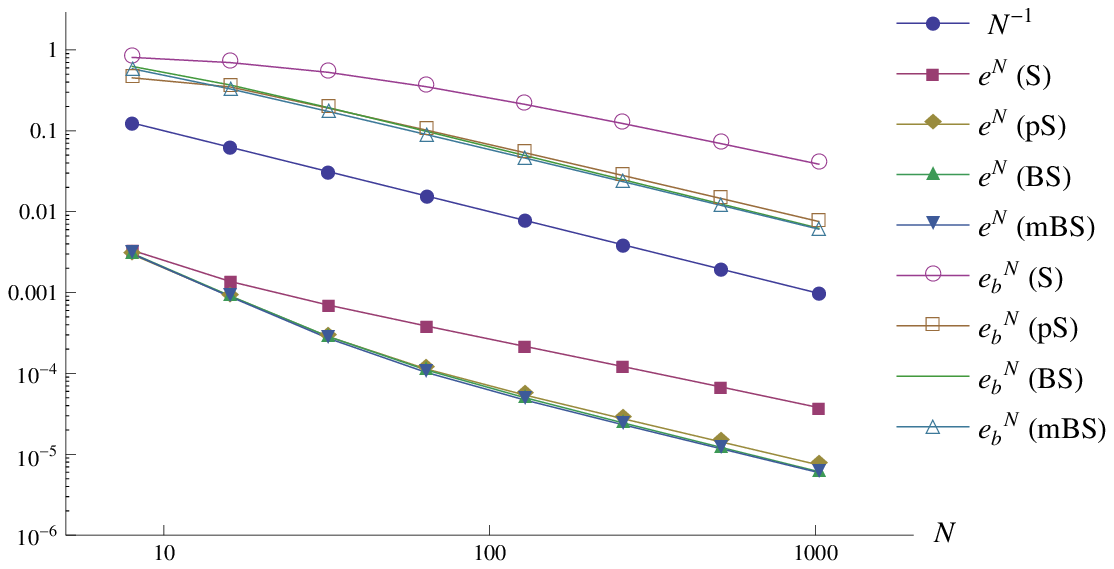}
\end{center}
\caption{Errors $e^N$ and $e^N_b$ on meshes from Table~\ref{tbl:mesh}, $\eps=2^{-20}$, $k=1$.}
\label{fig:k1}
\end{figure}

\begin{figure}[h!]
\begin{center}
\includegraphics[width=0.95\textwidth]{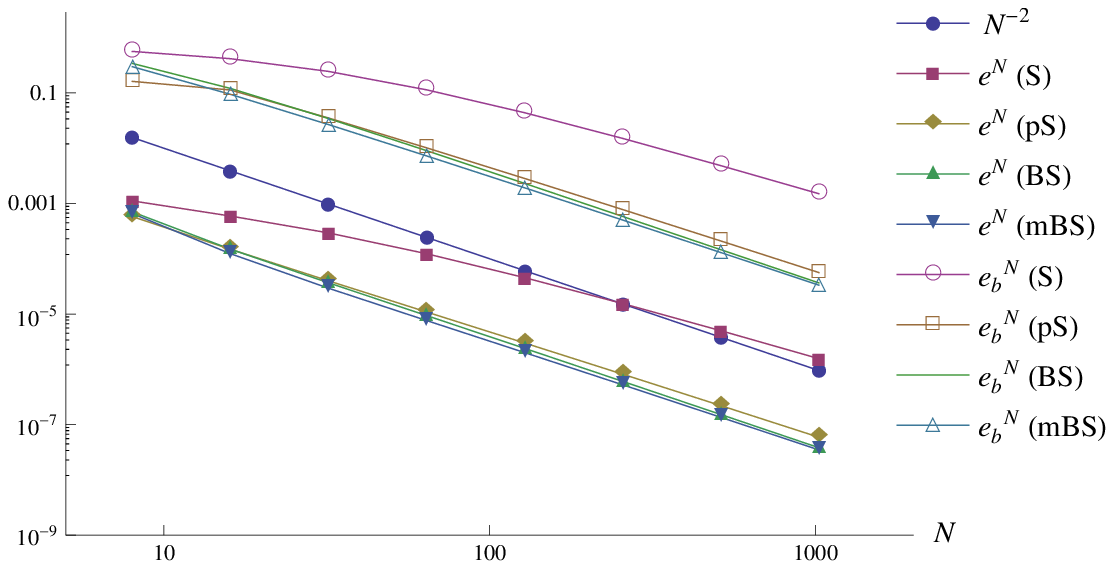}
\end{center}
\caption{Errors $e^N$ and $e^N_b$ on meshes from Table~\ref{tbl:mesh}, $\eps=2^{-20}$, $k=2$.}
\label{fig:k2}
\end{figure}

\begin{figure}[h!]
\begin{center}
\includegraphics[width=0.95\textwidth]{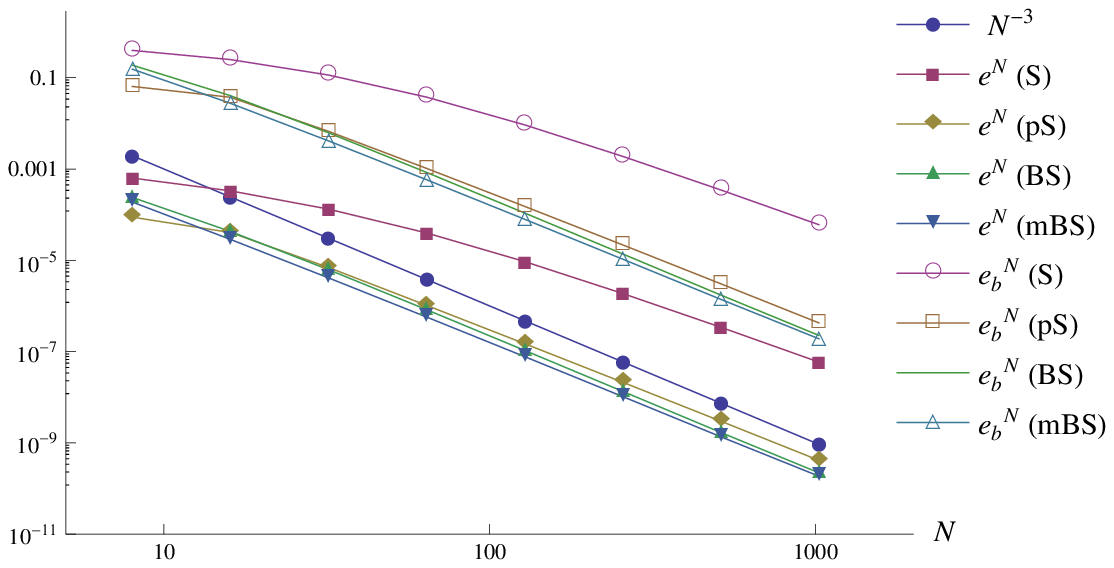}
\end{center}
\caption{Errors $e^N$ and $e^N_b$ on meshes from Table~\ref{tbl:mesh}, $\eps=2^{-20}$, $k=3$.}
\label{fig:k3}
\end{figure}

Let us first consider the meshes of Shishkin--type \eqref{eq:mesh}, with the mesh characterizing functions from Table~\ref{tbl:mesh}. For different values of $N$ and polynomial degrees $k=1,2,3$, in Table~\ref{tbl:k123} and Table~\ref{tbl:k123-B} we present the errors in the energy norm $e^N:=\|u-u^N\|_{dG}$ and in the balanced norm $e^N_b:=\|u-u^N\|_{dG,b}$, together with the rates of convergence estimated with the standard formulae
\begin{equation}\label{eq:rate}
p^N=\frac{\ln (e^N/e^{2N})}{\ln 2}, \qquad p^N_b=\frac{\ln (e_b^N/e_b^{2N})}{\ln 2}.\vspace{2mm}
\end{equation}
Except for the polynomial Shishkin mesh with $m>1$, all other meshes satisfy the assumption \eqref{eq:mesh_assump} that allows to bound the mesh width in the layer regions from below. Never\-theless, the numerical results confirm the order of convergence as proved in Theorem~\ref{thm:main} and Corollary~\ref{cor:B}. Comparison of the errors $e^N$ and $e^N_b$ on selected meshes is depicted on Figures~\ref{fig:k1}--\ref{fig:k3}.

\pagebreak

\begin{table}[h!]
\begin{center}
\caption{Energy norm error and rate of convergence, $\eps=2^{-20}$, $k=1,2,3$, for Shishkin (S--), polynomial Shishkin (pS--), Bakhvalov--Shishkin (BS--) and modified Bakhvalov--Shishkin (mBS--) mesh.\vspace{1mm}}
{\footnotesize
{\renewcommand{\arraystretch}{1.2}
\begin{tabular}{|c|cc|cc|cc|cc|} \hline
$N$ & \multicolumn{2}{c|}{S--mesh} & \multicolumn{2}{c|}{pS--mesh ($m=3$)} 
 	& \multicolumn{2}{c|}{BS--mesh} & \multicolumn{2}{c|}{mBS--mesh} \\ \hline
$k=1$ & $\|u-u^N\|_{dG}$ & $p^N$ & $\|u-u^N\|_{dG}$ & $p^N$ 
	& $\|u-u^N\|_{dG}$ & $p^N$ & $\|u-u^N\|_{dG}$ & $p^N$\\ \hline
%$2^3$ & $3.323(-3)$ & $1.279$ & $2.951(-3)$ & $1.716$ & $3.042(-3)$ & $1.745$ & $3.008(-3)$ & $1.765$ \\	
$2^4$ & $1.369(-3)$ & $0.970$ & $8.984(-4)$ & $1.660$ & $9.072(-4)$ & $1.670$ & $8.850(-4)$ & $1.707$ \\	
$2^5$ & $6.991(-4)$ & $0.859$ & $2.843(-4)$ & $1.327$ & $2.850(-4)$ & $1.372$ & $2.710(-4)$ & $1.400$ \\	
$2^6$ & $3.855(-4)$ & $0.832$ & $1.133(-4)$ & $1.065$ & $1.101(-4)$ & $1.131$ & $1.026(-4)$ & $1.126$ \\
$2^7$ & $2.165(-4)$ & $0.829$ & $5.415(-5)$ & $0.972$ & $5.027(-5)$ & $1.034$ & $4.702(-5)$ & $1.015$ \\
$2^8$ & $1.219(-4)$ & $0.838$ & $2.761(-5)$ & $0.952$ & $2.454(-5)$ & $1.007$ & $2.327(-5)$ & $0.987$ \\
$2^9$ & $6.818(-5)$ & $0.851$ & $1.427(-5)$ & $0.951$ & $1.221(-5)$ & $1.001$ & $1.174(-5)$ & $0.983$ \\
$2^{10}$ & $3.781(-5)$ & $-$ & $7.382(-6)$ & $-$ & $6.099(-6)$ & $-$	& $5.941(-6)$ & $-$ \\ \hline
$k=2$ & $\|u-u^N\|_{dG}$ & $p^N$ & $\|u-u^N\|_{dG}$ & $p^N$ 
	& $\|u-u^N\|_{dG}$ & $p^N$ & $\|u-u^N\|_{dG}$ & $p^N$\\ \hline
%$2^3$ & $1.116(-3)$ & $0.910$ & $5.802(-4)$ & $1.962$ & $7.045(-4)$ & $2.239$ & $6.588(-4)$ & $2.412$ \\	
$2^4$ & $5.940(-4)$ & $1.014$ & $1.489(-4)$ & $1.899$ & $1.493(-4)$ & $2.020$ & $1.237(-4)$ & $2.075$ \\	
$2^5$ & $2.942(-4)$ & $1.251$ & $3.994(-5)$ & $1.871$ & $3.680(-5)$ & $1.968$ & $2.937(-5)$ & $1.942$ \\	
$2^6$ & $1.236(-4)$ & $1.438$ & $1.092(-5)$ & $1.864$ & $9.407(-6)$ & $1.972$ & $7.647(-6)$ & $1.926$ \\
$2^7$ & $4.562(-5)$ & $1.557$ & $3.000(-6)$ & $1.875$ & $2.398(-6)$ & $1.983$ & $2.012(-6)$ & $1.938$ \\
$2^8$ & $1.550(-5)$ & $1.630$ & $8.179(-7)$ & $1.888$ & $6.065(-7)$ & $1.991$ & $5.252(-7)$ & $1.950$ \\
$2^9$ & $5.008(-6)$ & $1.679$ & $2.210(-7)$ & $1.899$ & $1.526(-7)$ & $1.995$ & $1.359(-7)$ & $1.960$ \\
$2^{10}$ & $1.564(-6)$ & $-$ & $5.927(-8)$ & $-$ & $3.827(-8)$ & $-$ & $3.494(-8)$ & $-$ \\ \hline
$k=3$ & $\|u-u^N\|_{dG}$ & $p^N$ & $\|u-u^N\|_{dG}$ & $p^N$ 
	& $\|u-u^N\|_{dG}$ & $p^N$ & $\|u-u^N\|_{dG}$ & $p^N$\\ \hline
%$2^3$ & $6.335(-4)$ & $0.928$ & $8.827(-5)$ & $1.117$ & $2.387(-4)$ & $2.475$ & $1.888(-4)$ & $2.715$ \\	
$2^4$ & $3.331(-4)$ & $1.330$ & $4.069(-5)$ & $2.546$ & $4.294(-5)$ & $2.781$ & $2.874(-5)$ & $2.800$ \\	
$2^5$ & $1.325(-4)$ & $1.734$ & $6.967(-6)$ & $2.756$ & $6.246(-6)$ & $2.913$ & $4.125(-6)$ & $2.848$ \\	
$2^6$ & $3.982(-5)$ & $2.085$ & $1.031(-6)$ & $2.786$ & $8.294(-7)$ & $2.962$ & $5.729(-7)$ & $2.882$ \\
$2^7$ & $9.387(-6)$ & $2.322$ & $1.495(-7)$ & $2.810$ & $1.065(-7)$ & $2.982$ & $7.771(-8)$ & $2.907$ \\
$2^8$ & $1.877(-6)$ & $2.457$ & $2.132(-8)$ & $2.831$ & $1.347(-8)$ & $2.989$ & $1.036(-8)$ & $2.924$ \\
$2^9$ & $3.419(-7)$ & $2.534$ & $2.996(-9)$ & $2.848$ & $1.697(-9)$ & $2.937$ & $1.365(-9)$ & $2.867$ \\
$2^{10}$ & $5.905(-8)$ & $-$ & $4.160(-10)$ & $-$ & $2.216(-10)$ & $-$ & $1.871(-10)$ & $-$ \\ \hline
\end{tabular} \label{tbl:k123}
}}
\end{center}
\end{table}

\pagebreak

\begin{table}[h!]
\begin{center}
\caption{Balanced norm error and rate of convergence, $\eps=2^{-20}$, $k=1,2,3$, for Shishkin (S--), polynomial Shishkin (pS--), Bakhvalov--Shishkin (BS--) and modified Bakhvalov--Shishkin (mBS--) mesh.\vspace{1mm}}
{\footnotesize
{\renewcommand{\arraystretch}{1.2}
\begin{tabular}{|c|cc|cc|cc|cc|} \hline
$N$ & \multicolumn{2}{c|}{S--mesh} & \multicolumn{2}{c|}{pS--mesh ($m=3$)} 
 	& \multicolumn{2}{c|}{BS--mesh} & \multicolumn{2}{c|}{mBS--mesh} \\ \hline
$k=1$ & $\|u-u^N\|_{dG,b}$ & $p^N_b$ & $\|u-u^N\|_{dG,b}$ & $p^N_b$ 
	& $\|u-u^N\|_{dG,b}$ & $p^N_b$ & $\|u-u^N\|_{dG,b}$ & $p^N_b$\\ \hline
%$2^3$ & $8.078(-1)$ & $0.213$ & $4.514(-1)$ & $0.374$ & $6.245(-1)$ & $0.763$ & $5.840(-1)$ & $0.828$ \\	
$2^4$ & $6.967(-1)$ & $0.400$ & $3.483(-1)$ & $0.871$ & $3.680(-1)$ & $0.931$ & $3.290(-1)$ & $0.926$ \\	
$2^5$ & $5.279(-1)$ & $0.592$ & $1.905(-1)$ & $0.904$ & $1.931(-1)$ & $0.973$ & $1.732(-1)$ & $0.948$ \\	
$2^6$ & $3.502(-1)$ & $0.721$ & $1.018(-1)$ & $0.924$ & $9.834(-2)$ & $0.988$ & $8.980(-2)$ & $0.960$ \\
$2^7$ & $2.125(-1)$ & $0.788$ & $5.364(-2)$ & $0.935$ & $4.959(-2)$ & $0.994$ & $4.616(-2)$ & $0.969$ \\
$2^8$ & $1.230(-1)$ & $0.824$ & $2.805(-2)$ & $0.943$ & $2.489(-2)$ & $0.997$ & $2.358(-2)$ & $0.975$ \\
$2^9$ & $6.950(-2)$ & $0.846$ & $1.459(-2)$ & $0.949$ & $1.247(-2)$ & $0.999$ & $1.199(-2)$ & $0.980$ \\
$2^{10}$ & $3.866(-2)$ & $-$ & $7.556(-3)$ & $-$ & $6.242(-3)$ & $-$ & $6.080(-3)$ & $-$ \\ \hline
$k=2$ & $\|u-u^N\|_{dG,b}$ & $p^N_b$ & $\|u-u^N\|_{dG,b}$ & $p^N_b$ 
	& $\|u-u^N\|_{dG,b}$ & $p^N_b$ & $\|u-u^N\|_{dG,b}$ & $p^N_b$\\ \hline
%$2^3$ & $5.644(-1)$ & $0.437$ & $1.620(-1)$ & $0.517$ & $3.425(-1)$ & $1.503$ & $2.986(-1)$ & $1.660$ \\	
$2^4$ & $4.168(-1)$ & $0.764$ & $1.132(-1)$ & $1.683$ & $1.209(-1)$ & $1.811$ & $9.447(-2)$ & $1.817$ \\	
$2^5$ & $2.454(-1)$ & $1.106$ & $3.526(-2)$ & $1.793$ & $3.444(-2)$ & $1.927$ & $2.681(-2)$ & $1.883$ \\	
$2^6$ & $1.140(-1)$ & $1.376$ & $1.018(-2)$ & $1.842$ & $9.055(-3)$ & $1.970$ & $7.269(-3)$ & $1.917$ \\
$2^7$ & $4.391(-2)$ & $1.543$ & $2.840(-3)$ & $1.869$ & $2.311(-3)$ & $1.987$ & $1.925(-3)$ & $1.937$ \\
$2^8$ & $1.506(-2)$ & $1.635$ & $7.775(-4)$ & $1.886$ & $5.829(-4)$ & $1.994$ & $5.025(-4)$ & $1.951$ \\
$2^9$ & $4.850(-3)$ & $1.688$ & $2.104(-4)$ & $1.898$ & $1.463(-4)$ & $1.997$ & $1.300(-4)$ & $1.960$ \\
$2^{10}$ & $1.505(-3)$ & $-$ & $5.643(-5)$ & $-$ & $3.666(-5)$ & $-$ & $3.342(-5)$ & $-$ \\ \hline
$k=3$ & $\|u-u^N\|_{dG,b}$ & $p^N_b$ & $\|u-u^N\|_{dG,b}$ & $p^N_b$ 
	& $\|u-u^N\|_{dG,b}$ & $p^N_b$ & $\|u-u^N\|_{dG,b}$ & $p^N_b$\\ \hline
%$2^3$ & $3.973(-1)$ & $0.669$ & $6.447(-2)$ & $0.790$ & $1.884(-1)$ & $2.219$ & $1.542(-1)$ & $2.485$ \\	
$2^4$ & $2.499(-1)$ & $1.124$ & $3.728(-2)$ & $2.478$ & $4.047(-2)$ & $2.699$ & $2.753(-2)$ & $2.734$ \\	
$2^5$ & $1.146(-1)$ & $1.608$ & $6.693(-3)$ & $2.695$ & $6.231(-3)$ & $2.886$ & $4.137(-3)$ & $2.827$ \\	
$2^6$ & $3.761(-2)$ & $2.019$ & $1.033(-3)$ & $2.764$ & $8.431(-4)$ & $2.954$ & $5.832(-4)$ & $2.876$ \\
$2^7$ & $9.276(-3)$ & $2.291$ & $1.521(-4)$ & $2.803$ & $1.088(-4)$ & $2.980$ & $7.944(-5)$ & $2.906$ \\
$2^8$ & $1.895(-3)$ & $2.444$ & $2.179(-5)$ & $2.829$ & $1.379(-5)$ & $2.988$ & $1.060(-5)$ & $2.923$ \\
$2^9$ & $3.484(-4)$ & $2.529$ & $3.066(-6)$ & $2.848$ & $1.737(-6)$ & $2.937$ & $1.398(-6)$ & $2.868$ \\
$2^{10}$ & $6.037(-5)$ & $-$ & $4.260(-7)$ & $-$ & $2.268(-7)$ & $-$ & $1.915(-7)$ & $-$ \\ \hline
\end{tabular} \label{tbl:k123-B}
}}
\end{center}
\end{table}

\pagebreak

\pagebreak

\vspace{5mm}

In order to show the robustness of the NIPG method on Shishkin--type meshes, we fix the number of degrees of freedom and the polynomial degree, and measure errors for various values of the perturbation parameter. Table~\ref{tbl:k2} shows the results only for Shishkin and Bakhvalov--Shishkin mesh with $N=2^{10}$ and $k=2$, with similar behavior on other S--type meshes. We observe the energy norm error $\|u-u^N\|_{dG}$ decreases when $\eps\to 0$, unlike the error in the balanced norm $\|u-u^N\|_{dG,b}$ that remains constant.

\begin{figure}[h!]
\begin{center}
\includegraphics[trim={0 0 4cm 0},clip,width=0.65\textwidth]{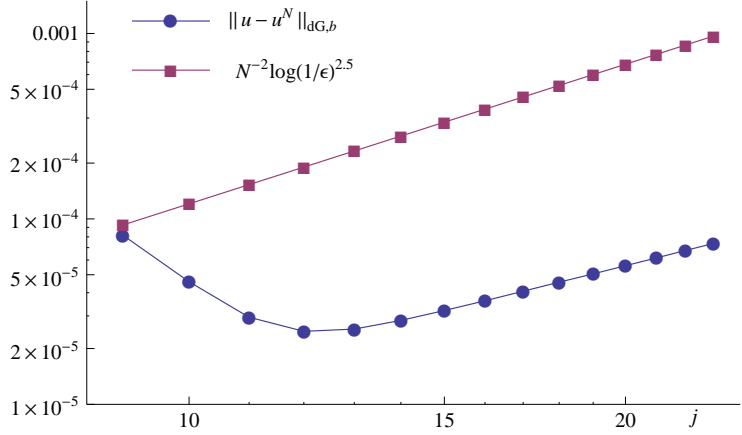}
\end{center}
\caption{Curve slopes for the balanced norm error $\|u-u^N\|_{dG,b}$ on DL--mesh and $N^{-k}\ln(1/\eps)^{k+1/2}$, for $N=N_{DL}=1024$, $k=2$ and $\eps=2^{-j}$.}\vspace{3mm}\label{fig:log}
\end{figure}

The results on the Duran--Lombardi mesh \eqref{eq:graded} are presented in Table~\ref{tbl:k123-G}. Here we choose $\eps=2^{-20}$, $k=1,2,3$, and $H=2^{-1},\dots,2^{-6}$. The number of mesh subintervals is denoted with $N_{DL}$. The rate of convergence in this case is estimated with $r^H$ and $r^H_b$ that are evaluated similarly to \eqref{eq:rate}, now taking the errors with the mesh parameters $H$ and $H/2$. Figure~\ref{fig:log} illustrates the influence of the logarithmic factor in error estimates on graded meshes. The results refer to balanced norm errors with $k=2$, $N_{DL}=1024$ and $\eps=2^{-9},\dots,2^{-23}$, where we notice similar slopes of the curves for $\|u-u^N\|_{dG,b}$ and $N^{-k}\ln(1/\eps)^{k+1/2}$, cf. Corollary~\ref{cor:2}.

\pagebreak

\begin{table}[h!]
\begin{center}
\caption{Errors $\|u-u^N\|_{dG}$ and $\|u-u^N\|_{dG,b}$, $N=2^{10}$, $k=2$, for Shishkin (S--) and Bakhvalov--Shishkin (BS--) mesh.\vspace{3mm}}
{\small
{\renewcommand{\arraystretch}{1.2}
\begin{tabular}{|c|cc|cc|} \hline
 & \multicolumn{2}{c|}{S--mesh} & \multicolumn{2}{c|}{BS--mesh} \\ \hline
$\eps$ & $\|u-u^N\|_{dG}$ & $\|u-u^N\|_{dG,b}$ & $\|u-u^N\|_{dG}$ & $\|u-u^N\|_{dG,b}$ \\ \hline
$2^{-10}$ & $5.048(-5)$ & $1.520(-3)$ & $2.263(-6)$ & $7.077(-5)$ \\
$2^{-11}$ & $3.549(-5)$ & $1.511(-3)$ & $1.101(-6)$ & $4.773(-5)$ \\
$2^{-12}$ & $2.505(-5)$ & $1.507(-3)$ & $6.594(-7)$ & $3.983(-5)$ \\
$2^{-13}$ & $1.770(-5)$ & $1.506(-3)$ & $4.422(-7)$ & $3.754(-5)$ \\
$2^{-14}$ & $1.251(-5)$ & $1.506(-3)$ & $3.080(-7)$ & $3.691(-5)$ \\
$2^{-15}$ & $8.846(-6)$ & $1.505(-3)$ & $2.169(-7)$ & $3.674(-5)$ \\
$2^{-16}$ & $6.254(-6)$ & $1.505(-3)$ & $1.532(-7)$ & $3.669(-5)$ \\
$2^{-17}$ & $4.422(-6)$ & $1.505(-3)$ & $1.083(-7)$ & $3.667(-5)$ \\
$2^{-18}$ & $3.127(-6)$ & $1.505(-3)$ & $7.654(-8)$ & $3.666(-5)$ \\
$2^{-19}$ & $2.211(-6)$ & $1.505(-3)$ & $5.412(-8)$ & $3.666(-5)$ \\
$2^{-20}$ & $1.564(-6)$ & $1.505(-3)$ & $3.827(-8)$ & $3.666(-5)$ \\ \hline
\end{tabular}\label{tbl:k2}
}}
\end{center}
\end{table}

\pagebreak

\begin{table}[h!]
\begin{center}
\caption{Errors $\|u-u^N\|_{dG}$ and $\|u-u^N\|_{dG,b}$, and rates of convergence, $\eps=2^{-20}$, $k=1,2,3$, for the Duran--Lombardi mesh.\vspace{3mm}}
{\small
{\renewcommand{\arraystretch}{1.2}
\begin{tabular}{|cc|cccc|} \hline
$H$ & $N_{DL}$ & \multicolumn{4}{c|}{DL--mesh} \\ \hline
\multicolumn{2}{|c|}{$k=1$} & $\|u-u^N\|_{dG}$ & $r^H$ & $\|u-u^N\|_{dG,b}$ & $r^H_b$ \\ \hline
$2^{-1}$ & $70$ & $4.933(-4)$ & $0.711$ & $1.505(-1)$ & $0.969$ \\
$2^{-2}$ & $128$ & $3.014(-4)$ & $1.887$ & $7.688(-2)$ & $0.981$ \\
$2^{-3}$ & $240$ & $8.146(-5)$ & $1.502$ & $3.895(-2)$ & $0.989$ \\
$2^{-4}$ & $468$ & $2.876(-5)$ & $1.359$ & $1.962(-2)$ & $0.994$ \\
$2^{-5}$ & $920$ & $1.121(-5)$ & $1.156$ & $9.852(-3)$ & $0.997$ \\
$2^{-6}$ & $1828$ & $5.032(-6)$ & $-$ & $4.937(-3)$ & $-$ \\ \hline
\multicolumn{2}{|c|}{$k=2$} & $\|u-u^N\|_{dG}$ & $r^H$ & $\|u-u^N\|_{dG,b}$ & $r^H_b$ \\ \hline
$2^{-1}$ & $70$ & $3.352(-5)$ & $2.191$ & $1.392(-2)$ & $1.815$ \\
$2^{-2}$ & $128$ & $7.340(-6)$ & $2.243$ & $3.956(-3)$ & $1.905$ \\
$2^{-3}$ & $240$ & $1.550(-6)$ & $2.085$ & $1.056(-3)$ & $1.953$ \\
$2^{-4}$ & $468$ & $3.653(-7)$ & $2.016$ & $2.728(-4)$ & $1.976$ \\
$2^{-5}$ & $920$ & $9.030(-8)$ & $1.996$ & $6.932(-5)$ & $1.988$ \\
$2^{-6}$ & $1828$ & $2.263(-8)$ & $-$ & $1.748(-5)$ & $-$ \\ \hline
\multicolumn{2}{|c|}{$k=3$} & $\|u-u^N\|_{dG}$ & $r^H$ & $\|u-u^N\|_{dG,b}$ & $r^H_b$ \\ \hline
$2^{-1}$ & $70$ & $1.689(-6)$ & $2.107$ & $7.900(-4)$ & $2.652$ \\
$2^{-2}$ & $128$ & $3.920(-7)$ & $3.767$ & $1.257(-4)$ & $2.792$ \\
$2^{-3}$ & $240$ & $2.879(-8)$ & $3.183$ & $1.815(-5)$ & $2.885$ \\
$2^{-4}$ & $468$ & $3.171(-9)$ & $3.183$ & $2.456(-6)$ & $2.918$ \\
$2^{-5}$ & $920$ & $3.491(-10)$ & $2.612$ & $3.250(-7)$ & $2.495$ \\
$2^{-6}$ & $1828$ & $5.710(-11)$ & $-$ & $5.764(-8)$ & $-$ \\ \hline
\end{tabular}\label{tbl:k123-G}
}}
\end{center}
\end{table}

\pagebreak

\section{Remarks to systems of reaction--diffusion\\ equations and the two--dimensional case}

A system of $s$ reaction--diffusion problems can be written in the form
\begin{equation}\label{eq:sys}
\left\{
\begin{array}{c}
-\eps^2\boldsymbol{u}''+\boldsymbol{A}\boldsymbol{u}
	=\boldsymbol{f}, \quad \text{in } (0,1), \\[1ex]
\boldsymbol{u}(0)=\boldsymbol{u}(1)=0,
\end{array}
\right.
\end{equation}
where the coupling matrix $\boldsymbol{A}:[0,1]\to\mathbb{R}^{s,s}$ is matrix--valued function and $\boldsymbol{u},\boldsymbol{f}:[0,1]\to\mathbb{R}^s$ are vector--valued. The Galerkin finite element method for the discretization of \eqref{eq:sys} with $s=2$ was first considered in \cite{LM04}, while a more general theory was devised in \cite{L09}. 

In a balanced norm, so far there exists only a result of Lin and Stynes \cite{LS15}. Following the basic idea from \cite{LS12}, but using $C^1-$elements instead of mixed finite elements, they introduce the bilinear form
\[
\eps^2(\boldsymbol{w}',\boldsymbol{v}')+(\boldsymbol{Aw},\boldsymbol{v})
	+\eps^3(\boldsymbol{w}'',\boldsymbol{v}'')+\eps((\boldsymbol{Aw})',\boldsymbol{v}')
\]
and analyse the finite element method for quadratic $C^1-$elements. The analysis for the Galerkin method with $C^0-$elements is open \cite{R17}.

For the discontinuous Galerkin method, however, our results can be extended to the system of reaction--diffusion equations \eqref{eq:sys}, as well as to the two--dimensional reaction--diffusion problem, \cite{RZ03}. It is more or less a technical question to generalize the results from \cite{RZ03}  to more general meshes, including error estimates in a balanced norm.

\section{Acknowledgements}

The research of the first author was supported by the Ministry of Education, Science and Technological Development of the Republic of Serbia, under Grant No.~174030.

%\bigskip
%\noindent {\bf References}


\begin{thebibliography}{99}

%{\small{

\bibitem{BZ16}
M.~Brdar, H.~Zarin, On graded meshes for a two--parameter singularly perturbed problem, Appl. Math. Comput. 282 (2016), 97--107.

\bibitem{CX15} 
P.~Constantinou, C.~Xenophontos, Finite element analysis of an exponentially graded mesh for singularly perturbed problems, CMAM 15(2) (2015), 135--143.

\bibitem{DL06} 
R.G.~Duran, A.I.~Lombardi, Finite element approximation of convection-diffusion problems on graded meshes, Appl. Num. Math. 56 (2006), 1314--1325.

\bibitem{FR13}
S.~Franz, H.--G.~Roos, Error estimation in a balanced norm for a convection--diffusion problem with two different boundary layers, Calcolo 51(3) (2014), 423--440.

\bibitem{FX16} 
S.~Franz, C.~Xenophontos, On a connection between layer-adapted exponentially graded and S--type meshes, submitted (2016). http://arxiv.org/abs/1611.07213

\bibitem{HSS02}
P.~Houston, Ch.~Schwab, E.~S\"uli, Discontinuous $hp-$finite element methods for advection--diffusion--reaction problems, SIAM J. Numer. Anal. 39(6) (2002), 2133--2163.

\bibitem{LS12}
R.~Lin, M.~Stynes, A balanced finite element method for singularly perturbed reaction--diffusion problems, SIAM J. Numer. Anal. 50(5) (2012), 2729--2743.

\bibitem{LS15} 
R.~Lin, M.~Stynes, A balanced finite element method for a system of singularly perturbed reaction--diffusion two--point boundary value problems, Numer. Algor. 70 (2015), 691--707.

\bibitem{L10}
T.~Linss, Layer--Adapted Meshes for Reaction--Convection--Diffusion Problems, Springer, 2010.

\bibitem{L09}
T.~Linss, Analysis of a FEM for a coupled system of singularly perturbed reaction--diffusion equations, Numer. Algor. 50 (2009), 283--291.

\bibitem{LM04}
T.~Linss, N.~Madden, A finite element analysis of a coupled system of singularly perturbed reaction--diffusion equations, Appl. Math. Comput. 148 (2004), 869--880.

\bibitem{M02}
J.M.~Melenk, $hp-$Finite Element Methods for Singular Perturbations, Vol. 1796 of Lecture Notes in Mathematics, Springer, Berlin, 2002.

\bibitem{MX16}
J.M.~Melenk, C.~Xenophontos, Robust exponential convergence of $hp-$FEM in balanced norms for singularly perturbed reaction--diffusion equations, Calcolo 53 (2016), 105--132.

\bibitem{MOS96}
J.J.H.~Miller, E.~OÕRiordan, G.~Shishkin, Fitted Numerical Methods for Singular Perturbation Problems, World Scientific, Singapore, 1996.

\bibitem{R17} 
H.--G.~Roos, Error estimates in balanced norms of finite element methods on layer-adapted meshes for
second order reaction--diffusion problems. Proc. of BAIL 2016, Beijing (to appear)

\bibitem{R12}
H.--G.~Roos, Robust numerical methods for singularly perturbed differential equations: a survey 2008-2012, International Scholarly Research Network, ISRN Applied Mathematics, Vol. 2012, Article ID 379547, 30 pages (doi:10.5402/2012/379547)

\bibitem{RL99}
H.--G.~Roos, T.~Linss, Sufficient conditions for uniform convergence on layer--adapted grids, Computing 63 (1999), 27--45.

\bibitem{RS14}
H.--G.~Roos, M.~Schopf, Convergence and stability in balanced norms of finite element methods on Shishkin meshes for reaction--diffusion problems, ZAMM 95 (2015), 551--565.

\bibitem{RST08}
H.--G.~Roos, M.~Stynes, L.~Tobiska, Robust Numerical Methods for Singularly Perturbed Differential Equations. Convection--Diffusion--Reaction and Flow Problems, Springer Series in Computational Mathematics, Springer--Verlag Berlin Heidelberg, 2008.

\bibitem{RZ03}
H.--G.~Roos, H.~Zarin, The Discontinuous Galerkin Finite Element Method for Singularly Perturbed Problems, in: Challenges in Scientific Computing -- CISC 2002. Proceedings of the Conference ÒChallenges in Scientific ComputingÓ, Berlin, Germany, October 2--5, 2002, (Springer, Berlin, 2003), pp. 246--267.

\bibitem{S05}
M.~Stynes, Steady--state convection--diffusion problems, Acta Numerica 14 (2005), 445--508. 

%}}
\end{thebibliography}
\end{document}